\documentclass[11pt]{amsart}

\def\diag{\mathop{\rm diag}}
\def\vol{\mathop{\rm Vol}}
\def\Jac{\mathop{\rm Jac}}
\def\Hom{\mathop{\rm Hom}}
\def\pr{\mathop{\rm pr}}

\def\half{{\frac{1}{2}}}

\def\zx{Z_{X_0}}
\def\tr{\mathop{\rm trace}}
\def\Cc{{\rm C_c}}
\def\inv{{^{-1}}}

\def\SL{{\rm SL}}
\def\GL{{\rm GL}}
\def\M{{\rm M}}
\def\SO{{\rm SO}}
\def\O{{\rm O}}

\def\calF{\mathcal{F}}

\def\cF{\mathfrak{F}}

\def\cO{\mathfrak{O}}
\def\SM{\mathcal{S\!M}}
\def\SP{\mathcal{S\!P}}
\def\cP{\mathcal{P}}
\def\cM{\mathcal{M}}
\def\frakM{\mathfrak{M}}
\def\cU{\mathfrak{U}}

\def\bfa{\mbox{\boldmath{$a$}}}
\def\bfb{\mbox{\boldmath{$b$}}}
\def\bfc{\mbox{\boldmath{$c$}}}
\def\bfn{\mbox{\boldmath{$n$}}}

\def\bfs{\mbox{\boldmath{$s$}}}
\def\bft{\mbox{\boldmath{$t$}}}
\def\sbft{\mbox{\boldmath{$\scriptstyle t$}}}
\def\bfu{\mbox{\boldmath{$u$}}}

\def\bfx{\mbox{\boldmath{$x$}}}
\def\sbfx{\mbox{\boldmath{$\scriptstyle x$}}}
\def\bfy{\mbox{\boldmath{$y$}}}
\def\bfz{\mbox{\boldmath{$z$}}}

\def\bfomega{\mbox{\boldmath{$\omega$}}}

\def\N{{\mathbb{N}}}
\def\Z{{\mathbb{Z}}}
\def\Q{{\mathbb{Q}}}
\def\R{{\mathbb{R}}}
\def\C{{\mathbb{C}}}

\def\cl#1{{\overline{#1}}}

\def\con#1{{\,^{#1}\!}}
\def\trn#1{{\,^{\bf t}\!#1}}

\newtheorem{theo}{Theorem}[section]
\newtheorem{lema}[theo]{Lemma}
\newtheorem{prop}[theo]{Proposition}

\newtheorem{rem}{Remark}[section]
\newtheorem{notation}{Notation}[section]

\def\bp{\noindent{\it Proof. }}
\def\ep{\qed}

\begin{document}

\title[Counting integral matrices]{Counting integral matrices with a
given characteristic polynomial} 

\author{Nimish A. Shah}

\address{Tata Institute of Fundamental Research, Homi Bhabha Rd.,
Mumbai 400005, India; e-mail: nimish@math.tifr.res.in}
 
\date{\today}

\begin{abstract}
We give a simpler proof of an earlier result giving an asymptotic
estimate for the number of integral matrices, in large balls, with a
given monic integral irreducible polynomial as their common
characteristic polynomial. The proof uses equidistributions of
polynomial trajectories on $\SL(n,\R)/\SL(n,\Z)$, which is a generalization
of Ratner's theorem on equidistributions of unipotent trajectories.

We also compute the exact constants appearing in the above mentioned
asymptotic estimates.
\end{abstract}

\maketitle

\section{Introduction}
Let $P$ be a monic polynomial of degree $n$ ($n\geq 2$) with integral
coefficients which is irreducible over $\Q$. Let
\[
V_P=\{X\in\M_n(\R):\det(\lambda I-X)=P(\lambda)\}.
\] 
Since $P$ has $n$ distinct roots, $V_P$ is the set of real $n\times
n$-matrices $X$ such that roots of $P$ are the eigenvalues of $X$. Let
$V_P(\Z)$ denote that set of matrices in $V_P$ with integral entries.
Let $B_T$ denote the ball in $M_n(\R)$ centred at $0$ and of radius
$T$ with respect to the Euclidean norm: $\|(x_{ij})\|=(\sum_{i,j}
x_{ij}^2)^{1/2}$. We are interested in estimating, for large $T$, the
number of integer matrices in $B_T$ with characteristic polynomial
$P$.

\begin{theo}[\cite{EMS-ann}] \label{thm:EMS:counting}
There exists a constant $C_P>0$ such that 
\[
\lim_{T\to\infty} \frac{\#(V_P(\Z)\cap B_T)}{T^{n(n-1)/2}} = C_P.
\]
\end{theo}

A formula for $C_P$, in the general case, is given in
Theorem~\ref{thm:C_P}. Under an additional hypothesis, the formula for
$C_P$ is simpler and it can be given as follows (Cf.~\cite{EMS-ann}):

\begin{theo} \label{thm:Z[alpha]}
Let $\alpha$ be a root of $P$ and $K=\Q(\alpha)$. Suppose that
$\Z[\alpha]$ is the integral closure of $\Z$ in $K$. Then
\[
C_P=\frac{2^{r_1}(2\pi)^{r_2}hR}{w\sqrt{D}}\cdot
\frac{\pi^{m/2}/\Gamma(1+(m/2))}
{\prod_{s=2}^n\pi^{-s/2}\Gamma(s/2)\zeta(s)},
\]
where $h$ = ideal class number of $K$, $R$ = regulator of $K$, $w$ =
order of the group of roots of unity in $K$, $D$ = discriminant of
$K$, $r_1$ (resp.\ $r_2$) = number of real (resp.\ complex) places
of $K$, and $m=n(n-1)/2$. 
\end{theo}

\begin{rem} \label{rem:identify} The three components of the above
formula for $C_P$ are volumes of certain standard entities in geometry
of numbers (with respect to the canonical volume forms on the
respective spaces):
\begin{eqnarray*}
\vol(J^0(K)/K^\times) &=& \frac{2^{r_1}(2\pi)^{r_2}hR}{w\sqrt{D}},\\
\vol(B^m) &=& \pi^{m/2}/\Gamma(1+(m/2)), \\ 
\vol(\SM_n) &=& 
\prod_{s=2}^n \pi^{-s/2}\Gamma(s/2)\zeta(s).
\end{eqnarray*}
Here $J^0(K)/K^\times$ = the group of principal ideles of $K$ modulo
$K^\times$ (see~\cite[Sect.\ 5.4]{Koch}), $B^m$ = the unit ball in
$\R^m$, and $\SM_n$ = the determinant one surface in the Minkowski
fundamental domain $M_n$ in the space of $n\times n$ real positive
symmetric matrices with respect to the action of $\GL_n(\Z)$
(see~\cite[Sect.~4.4.4]{Terras}).
\end{rem}
  
\begin{rem} \label{rem:cyclotomic} 
The hypothesis of Theorem~\ref{thm:Z[alpha]} is satisfied if $\alpha$
is a root of unity (see \cite[Theorem~1.61]{Koch}).
\end{rem}

The conclusion of Theorem~\ref{thm:Z[alpha]} was obtained in
\cite{EMS-ann} under a further hypothesis that all roots of $P$ are
real.

In \cite{EMS-ann}, the proof of Theorem~\ref{thm:EMS:counting} is
based on the following: (1) the existence of limits of large
translates of certain algebraic measures as proved in \cite{EMS-gafa};
(2) showing that such limiting distributions are actually algebraic
measures, using Ratner's description of ergodic invariant measures of
unipotent flows \cite{Ratner-measure}; and (3) the verification that
certain condition, called the \emph{non-focusing condition}\/, holds in
the case of Theorem~\ref{thm:EMS:counting}. (See \cite{Ratner:ICM}).

A main purpose of this article is to provide a simple and a direct
proof of this theorem using the following result on equidistributions of
`polynomial like' trajectories on $\SL_n(\R)/\SL_n(\Z)$:

\begin{theo} \label{thm:polynomial}
Let $\Gamma$ be a lattice in $\SL_n(\R)$, $\mu$ the
$\SL_n(\R)$-invariant probability measure on $\SL_n(\R)/\Gamma$, and
$x\in\SL_n(\R)/\Gamma$. Let 
\[
\Theta=(\Theta_{ij})_{i,j=1}^n: \R^m\to \SL_n(\R)
\]
be a map such that each $\Theta_{ij}$ is a real valued polynomial in
$m$ variables, and $\Theta(0)=I$, the identity matrix. Suppose that
$\Theta(\R^m)$ is not contained in any proper closed subgroup $L$ of
$\SL_n(\R)$ such that the orbit $Lx$ is closed. Then for any $f\in
\Cc(\SL_n(\R)/\Gamma)$,
\[
\lim_{T\to\infty}\frac{1}{\vol(B(T))} \int_{B(T)}
f(\Theta(\bfs)x)\; d\bfs 
= \int f\; d\mu,
\]
where $B(T)$ denotes the ball of radius $T$ in $\R^m$ centered at $0$.

For $0\leq r\leq m$, let $B^+(T)=B(T)\,\cap\, (\R_+)^r\!\times\!\R^{m-r}$. Then
\[
\lim_{T\to\infty} \frac{1}{\vol(B^+_T)} \int_{B^+_T} 
f(\delta(\bfs)x)\; d\bfs 
= \int f\; d\mu, \qquad \forall f\in \Cc(\SL_n(\R)/\Gamma),
\]
where 
$\delta(\bfs):=\Theta(\sqrt{s_1},\ldots,\sqrt{s_r},s_{r+1},\ldots,s_m)$,
$\forall\, \bfs=(\R_+)^r\times R^{m-r}$.
\end{theo}

The first part of the theorem is a particular case of
\cite[Corollary~1.1]{Shah-polynomial}, whose proof can be readily
modified to prove the second part. This result is a generalization of
Ratner's theorem on equidistribution of orbits of one-dimensional
unipotent flows \cite{Ratner-distribution}. The main ingredient in its
proof is, just as in \cite{Ratner-distribution}, the classification of
ergodic invariant measures for unipotent flows.

Another purpose of this article is to obtain an expression for $C_P$
in terms of algebraic number theoretic constants associated with $P$;
this is carried out in Section~\ref{sec:volumes}.

As in \cite{EMS-ann}, the first step in the proof of
Theorem~\ref{thm:EMS:counting} is its reformulation to a question in
ergodic theory of subgroup actions on homogeneous spaces of Lie
groups; we follow the approach of Duke, Rudnick and Sarnak \cite{DRS}.

The second step is to reduce this question to one about
equidistribution of polynomial trajectories, so that
Theorem~\ref{thm:polynomial} can be applied. 

\section{Reduction to a question in ergodic theory}
\label{subsec:reduction} 
 
We write
\[
\con{g}X:=gXg\inv, \qquad  \forall g\in \GL_n(\R),\, \forall X\in
\M_n(\R).
\]
Put 
\[
\Gamma=\GL_n(\Z).
\]
If $X\in V_P(\Z)$ and $\gamma\in\Gamma$, then $\con{\gamma} X\in
V_P(\Z)$; and we denote the $\Gamma$-orbit through $X$ by
\[
\con{\Gamma}X:=\{\con{\gamma}X:\gamma\in\Gamma\}.
\]

Using a correspondence between $\Gamma$-orbits and ideal classes due
to Latimer and MacDuffee~\cite{MacDuffee}, in view of the finiteness of
class numbers of orders, one has the following: (see
Proposition~\ref{prop:finite:orbits2}).

\begin{prop}[Latimer and MacDuffee] \label{prop:finite:orbits}
There are only finitely many distinct $\Gamma$-orbits in $V_P(\Z)$.
\end{prop}

\begin{rem} \label{rem:B:H-C}
The above proposition is a particular case of a much general
`finiteness theorem' due to Borel and Harish-Chandra~\cite{B:H-C}.
\end{rem}

By Proposition~\ref{prop:finite:orbits}, to prove
Theorem~\ref{thm:EMS:counting} it is enough to prove the following.

\begin{theo} \label{thm:EMS:orbit}
Let $X\in V_P(\Z)$. Then there exists $c_X>0$ such that
\[
\lim_{T\to\infty} \frac{\#(\con{\Gamma}X\cap B_T)}{T^{n(n-1)/2}} = c_X.
\]
\end{theo}

\subsection{Considering a fixed  \(\Gamma\)-orbit.} 

Put $G=\{g\in \GL_n(\R): \det g=\pm 1\}$. Since the conjugation action
of $\GL_n(\R)$ on $V_P$ is transitive, the same holds for the action
of $G$ on $V_P$. Note that $\Gamma=\GL_n(\Z)$ is a lattice in $G$. Fix
any $X_0\in V_P(\Z)$. Put
\[
H=\{g\in G:\con{g}X_0=X_0\}.
\]
Then $H$ is a real algebraic torus defined over $\Q$. In
Section~\ref{sec:compact}, using the Dirichlet's unit theorem will
show the following.

\begin{prop} \label{prop:H:compact} 
$H/H\cap\Gamma$ is compact.
\end{prop}

\medskip
Define 
\[
R_T=\{g\in G: \con{g}X_0\in B_T\}/H \, \subset \, G/H,
\]
and $\chi_T$ denote its characteristic function. Then 
\begin{equation} \label{eq:R_T}
\#(\con{\Gamma}X_0\cap B_T) = \#(\Gamma[H] \cap R_T)
=\sum_{\dot{\gamma}\in\Gamma/\Gamma\cap H}\chi_T(\gamma[H]).
\end{equation}

We choose Haar measures $\tilde \mu$ (resp.\ $\tilde \nu$) on $G$
(resp.\ $H$). Let $\mu$ (resp.\ $\nu$) denote the left invariant
measure on $G/\Gamma$ (resp.\ $H/H\cap\Gamma$) corresponding to the
measure $\tilde\mu$ (resp.\ $\tilde\nu$).

Let $\eta$ be the corresponding left $G$-invariant measure on $G/H$
(see~\cite[Lemma~1.4]{Raghunathan}); that is, $\forall\, f\in\Cc(G)$,
\begin{equation} \label{eq:eta:def}
\int_G f\;d\tilde\mu=
\int_{gH\in G/H} \left(\int_H f(gh)\;
d\tilde\nu(h)\right)\;d\eta(gH).
\end{equation}

In Section~\ref{subsec:vol(R_T)} we show that there exists a constant
$c_\eta>0$ (see~\ref{eq:R_T:vol}) depending on $X_0$ such that
\begin{equation} \label{eq:vol(R_T)}
\lim_{T\to\infty} \eta(R_T)/T^{n(n-1)/2}=c_\eta.
\end{equation}

For all $T>0$ and $g\in G$, let
\begin{equation} \label{eq:F_T}
F_T(g\Gamma) := \#(g\Gamma[H]\cap R_T)
= \sum_{\dot{\gamma}\in \Gamma/(\Gamma\cap H)} \chi_T(g\gamma H).
\end{equation}
Note that $F_T$ is bounded, measurable, and vanishes outside a
compact set in $G/\Gamma$. By (\ref{eq:R_T}) and (\ref{eq:vol(R_T)}),
in order to prove Theorem~\ref{thm:EMS:orbit}, it is enough to prove
the following:

\begin{theo} \label{thm:EMS:R_T} \ 
\[
\lim_{T\to\infty} \frac{F_T(e\Gamma)}{\eta(R_T)} =
\frac{\nu(H/H\cap\Gamma)}{\mu(G/\Gamma)}.
\]
\end{theo}

From the computations in Section~\ref{sec:R1} and \ref{sec:D1}, one
can deduce the following: Given any $\kappa>1$ there exists a
neighbourhood $\Omega$ of $e$ in $G$ such that
\begin{equation} \label{eq:kappa:Omega}
R_{\kappa\inv T}\subset \Omega R_T \subset R_{\kappa T}.
\end{equation}

Now by (\ref{eq:vol(R_T)}), 
\begin{equation} \label{eq:kappa:T}
\lim_{\kappa\to 1}
\lim_{T\to\infty} \eta(R_{\kappa T})/\eta(R_T) =
\frac{\nu(H/H\cap\Gamma)}{\mu(G/\Gamma)}.
\end{equation}

By (\ref{eq:kappa:Omega}) and (\ref{eq:kappa:T}), 
in order to prove Theorem~\ref{thm:EMS:R_T}, it is enough to
prove the following weak convergence:

\begin{theo} \label{thm:F_T:weak} 
For any $f\in\Cc(G/\Gamma)$,
\[
\lim_{T\to\infty} \frac{\langle f,F_T\rangle}{\eta(R_T)} =
\frac{\nu(H/H\cap\Gamma)}{\mu(G/\Gamma)}\cdot \langle f,1\rangle.
\]
\end{theo}

Using Fubini's theorem we have the following:

\begin{prop}[\cite{DRS,EM}] \label{prop:fubini1}
For any $f\in\Cc(G/\Gamma)$,
\begin{eqnarray}
\langle f, F_T\rangle & = & \int_{G/\Gamma}
f(g\Gamma) \left(\sum_{\dot{\gamma}\in \Gamma/(H\cap\Gamma)}
\chi_T(g\gamma H)\right) \; d\mu(\dot{g}) \nonumber 
\\ & = &
\int_{G/H\cap\Gamma} f(g\Gamma) \chi_T(gH)
\;d\bar{\mu}(\dot{g}) \nonumber 
\\ & = &
\int_{R_T}
\left(\int_{H/H\cap\Gamma}
f(gh\Gamma)\; d\nu(\dot{h})\right) \,d\eta(\dot{g}),
\label{eq:fubini1}
\end{eqnarray}
where $\bar{\mu}$ is the left $G$-invariant measure on
$G/(H\cap\Gamma)$ corresponding to $\tilde\mu$, and $\dot{x}$ denotes
the appropriate coset of $x$.
\end{prop}

In \cite{EMS-ann} further analysis of the limit was carried out by
showing that, as $T\to\infty$, for `almost all' sequences
$g_iH\to\infty$ in $R_T$, the integral in the bracket of
Equation~\ref{eq:fubini1} converges to
$\frac{\nu(H/H\cap\Gamma)}{\mu(G/\Gamma)}\, \langle
f,1\rangle$. This then implies Theorem~\ref{thm:F_T:weak}.

In this article, our approach is to change the order of integration in
(\ref{eq:fubini1}), and then apply Theorem~\ref{thm:polynomial} to
find the limit. For this purpose, we need an explicit description of
$R_T$, and of the measure $\eta$.

\section{Integration on  \(R_T\)} \label{sec:2}

\begin{notation} \label{not:sigma}
Let $r_1$ be the number of real roots of $P$ and $r_2$ be the number
of pairs of complex conjugate roots of $P$. Since $P$ is irreducible,
all roots of $P$ are distinct, and $n=r_1+2r_2$. Fix a root $\alpha$
of $P$. Let $\sigma_i$ ($i=1,\ldots,r_1$) be
the distinct real embeddings of $\Q(\alpha)$. Let $\sigma_{r_1+i}$
($i=1,\ldots,2r_2$) be the distinct complex embeddings of
$\Q(\alpha)$, such that
\begin{equation} \label{eq:sigma:bar}
\sigma_{r_1+r_2+i}=\overline{\sigma_{r_1+i}}, \qquad 1\leq i\leq r_2.
\end{equation}
\end{notation}

Put 
\begin{equation} \label{eq:di}
d_i= \left\{\begin{array}{ll}  
\sigma_i(\alpha) & \mbox{ if $1\leq i\leq r_1$} \\
\left(
\begin{matrix}  
a_{i-r_1} & - b_{i-r_1} \\ 
b_{i-r_1} & a_{i-r_1} 
\end{matrix} 
\right) & \mbox{ if $r_1 < i \leq r_1+r_2$}, 
\end{array} \right. 
\end{equation}
where $a_i+b_i\sqrt{-1}:=\sigma_{r_1+i}(\alpha)$, $i=1,\ldots,r_2$.

\subsection{Diagonalization of  \(X\) and  \(H\).} 
\label{sec:diagonalize}

Let
\[
\begin{array}{lcl}  
X_1 &=& \diag(d_1,\ldots,d_{r_1+r_2})  \\[3pt] 
H_1 &=& \{g\in G:\con{g}X_1=X_1\}  \\[3pt] 
R^1_T &=& \{g\in G: \con{g}X_1 \in B_T\}/H_1.  
\end{array}
\]

Since the eigenvalues of $X_1$ are same as the roots of $P$, $X_1\in
V_P$. Let $g_0\in G$ be such that $\con{g_0}X_0=X_1$. 

Define $\psi:G\to G$ as $\psi(g)=g_0gg_0\inv$, $\forall\, g\in G$.
Then $H_1 = \psi(H)$ and $\psi_\ast(\tilde\mu)=\tilde\mu$.  We choose
a Haar measure $\tilde\nu_1$ on $H_1$ defined by
\begin{equation}\label{eq:tilde:nu1:def}
\tilde\nu_1:=\psi_\ast(\tilde\nu).
\end{equation}

Define $\bar\phi:G/H \to G/H_1$ as $\bar\phi(gH)=gg_0\inv H_1$,
$\forall\, g\in G$. Let $\eta_1:=\bar\phi_\ast(\eta)$. Then by
(\ref{eq:eta:def}), $\forall\, f\in\Cc(G)$,
\begin{equation} \label{eq:eta1:def}
\int_G f\; d\tilde\mu = \int_{G/H_1} \left(\int_{H_1}
f(gh_1)\;\tilde\nu_1(h_1)\right) \;d\eta_1(gH_1).
\end{equation}
Also 
\begin{equation} \label{eq:eta1:eta}
R^1_T=\bar\phi(R_T) \qquad \mbox{and} \qquad \eta_1(R^1_T)=\eta(R_T).
\end{equation}

Put $\Gamma_1=\psi(\Gamma)$. Define $\bar\psi:G/\Gamma\to G/\Gamma_1$
as $\bar\psi(g\Gamma)=\psi(g)\Gamma_1$, $\forall\, g\in G$.  Let
$\mu_1:=\bar\psi_\ast(\mu)$ and $\nu_1:=\bar\psi_\ast(\nu)$. Then
$\mu_1$ is the $G$-invariant measure on $G/\Gamma_1$ associated to
$\tilde\mu$. Also $\nu_1$ is the $H_1$-invariant measure on
\[
H_1/(H_1\cap \Gamma_1)\cong H_1\Gamma_1/\Gamma_1 =
\bar\psi(H\Gamma/\Gamma)
\]
associated to $\tilde \nu_1$, and
\begin{equation} \label{eq:nu1=nu}
\nu_1(H_1/H_1\cap\Gamma_1) = \nu(H/H\cap\Gamma).
\end{equation}

Now can rewrite Proposition~\ref{prop:fubini1} as follows:

\begin{prop} \label{prop:R1}
$\forall f\in\Cc(G/\Gamma)$, and
$f_1:=f\circ\bar{\psi}\inv\in\Cc(G/\Gamma_1)$, 
\begin{eqnarray*} 
\langle f,F_T\rangle   
&=& \int_{R_T}  
\left(\int_{H/H\cap\Gamma} f(gh\Gamma)\;d\nu(\dot{h})\right)
\;d\eta(\dot{g})\\
&=& \int_{R^1_T} \left(\int_{H_1/H_1\cap\Gamma_1}
f_1(gh\Gamma_1)\; d\nu_1(\dot{h})\right) \; d\eta_1(\dot{g}).
\end{eqnarray*}
\end{prop}

Due to this proposition, instead of integrating on $R_T$, it suffices
to integrate on $R^1_T$. Therefore we describe the measure $\eta_1$ on
$G/H_1$. For this purpose we want to express $G$ as $G=YH_1$, where
$Y$ is a product of certain subgroups and subsemigroups of $G$ (see
(\ref{eq:G=KN1UH1})). Later, in Section~\ref{sec:Haar} we will
decompose the Haar measure of $G$ into products of appropriate Haar
measures on these subgroups. This will allow us to describe $\eta_1$
as a product of the chosen Haar measures on the subgroups and
subsemigroups, whose product is $Y$
(Proposition~\ref{prop:G/H1:Haar}).

\subsection{Product decompositions of  \(G\).} \label{sec:decompose} 
In view of the above, first we will describe various subgroups of $G$,
and then obtain different product decompositions of $G$ into those
subgroups and their subsemigroups.

Let $\O(n)$ denote the group of orthogonal matrices in
$\GL_n(\R)$. Let
\begin{eqnarray}
N & = & \{\bfn:=(n_{ij})_{i,j=1}^n: n_{ij}\in\R,\, n_{ij}=0 \mbox{ if
$i>j$},\, n_{ii}=1\} \label{eq:N} \\
A & = &
\{\bfa:=\diag(a_1,\ldots,a_n):a_i>0,\,\prod_{i=1}^n a_i = 1\}.
\label{eq:A}
\end{eqnarray}

By Iwasawa decomposition, the map 
\[
(k,n,a)\mapsto kna \,:\, \O(n)\times N\times A\to G
\]
is a diffeomorphism.

For $i,j=1,\ldots,r_1+r_2$, let
\begin{equation} \label{eq:Mij:def}
M_{ij}=\scriptstyle{\left\{
\begin{array}{ll} 
\R & \mbox{if $i\leq r_1$, $j\leq r_1$}  \\
\M_{1\times 2}(\R) & \mbox{if $i\leq r_1$, $j>r_1$} \\ 
\M_{2\times 1}(\R) & \mbox{if $i>r_1$, $j\leq r_1$} \\ 
\M_2(\R) & \mbox{if $i>r_1$, $j>r_1$.}
\end{array} 
\right.}
\end{equation}
It will be convenient to express $g\in\M_n(\R)$ as
$g=(g_{ij})_{i,j=1}^{r_1+r_2}$, where $g_{ij}\in M_{ij}$. 

Put
\[
\begin{array}{cl}
\cU & =  \left(\prod_{1\leq i<j\leq r_1+r_2} M_{ij}\right) 
\cong \; \R^{\half n(n-1)-r_2}, \\[3pt]
u(\bfx)& =  (u_{ij});\ \bfx=(x_{ij}) \in\cU,\ 
M_{ij}\ni
u_{ij}=\scriptstyle{\left\{ \begin{array}{ll}  
0 & \mbox{ if $i>j$} \\ 
1 & \mbox{ if $i=j$} \\ 
x_{ij} & \mbox{ if $i<j$,} 
\end{array} 
\right.} \\
h(t) & =\left(\begin{matrix} 1 & t \cr 0 & 1
\end{matrix}\right),\ \forall\, t\in\R.
\end{array}
\]

Define
\begin{eqnarray} 
\qquad L_1 &=& \{\diag(1,\ldots,1,g_1,\ldots,g_{r_2})\in G:g_i\in\SL_2(\R)\} 
\nonumber\\
\qquad K_1 &=& \{\diag(1,\ldots,1,k_1,\ldots,k_{r_2})\in G:k_i\in\SO(2)\}
\nonumber \\[3pt]
\qquad N_1 &=& \{h(\bft)=\diag(1,\ldots,1,h(t_1),\ldots,h(t_{r_2}): 
\nonumber   \\
\ &\ & \qquad \ \qquad \bft=(t_i)\in\R^{r_2}\} \nonumber \\[3pt]
\qquad A_1 &=& \{\bfa_1=\diag(1,\ldots,1,b_1,\ldots,b_{r_2}): 
\label{eq:A1}\\
\ &\ & \qquad \ \qquad b_i=\diag(\beta_i,\beta_i\inv),\,\beta_i>0\} 
\nonumber \\[3pt]
\qquad U &=&
\{u(\bfx): \bfx=(x_{ij})\in\cU\}
\nonumber \\[3pt]
\qquad C &=& 
\{\bfc=
\diag(c_1,\ldots,c_{r_1},c_{r_1+1}^{1/2} I_2,\ldots,c_{r_1+r_2}^{1/2} I_2)
\in G: \label{eq:C} \\ 
\ &\ & \qquad \ \qquad c_i>0,\, \mbox{$\prod_{i=1}^{r_1+r_2} c_i=1$}\} 
\nonumber \\[3pt]
\qquad \Sigma &=&
\{\diag(\epsilon_1,\ldots,\epsilon_{r_1},I_2,\dots,I_2)\in G:
\epsilon_i=\pm 1\} 
\label{eq:Sigma}
\end{eqnarray}

We have the following product decompositions: 
\begin{equation} \label{eq:decompose}
\begin{array}{ll}
N=N_1\cdot U, & A=A_1\cdot C, \\[3pt]
 H_1=\Sigma \cdot K_1 \cdot C, &  
L=K_1\cdot N_1\cdot A_1.
\end{array}
\end{equation}
In each of the above decompositions, the product map, from the direct
product of the subgroups on the right hand side to the group on the
left hand side, is a diffeomorphism. We also note that
\begin{equation} \label{eq:N(U)}
\begin{array}{ll}
\Sigma\cdot C\subset Z_G(L), & N_G(U)=\Sigma \cdot C \cdot L\cdot U.
\end{array}
\end{equation}

Therefore
\begin{equation} \label{eq:G=KLUA}
\begin{array} {ll}
G & = \O(n)NA = \O(n)K_1\cdot N_1 U\cdot A_1 C \\[3pt]
\ & = \O(n)\cdot K_1N_1A_1\cdot UC \\[3pt]
\ & =  \O(n)\cdot L\cdot U\cdot C. 
\end{array}
\end{equation}

One has that $\SL_2(\R)=\SO(2)\cdot h(\R_+)\cdot\SO(2)$ (see
Proposition~\ref{prop:KNA=KNK}). Since $L\cong(\SL_2(\R))^{r_2}$, we
have that
\[
L=K_1 N_1^+K_1,
\]
where $N_1^+=\{h(\bft):\bft\in(\R_+)^{r_2}\}$. Now, in view of
(\ref{eq:decompose})--(\ref{eq:G=KLUA}), we have
\begin{equation} \label{eq:G=KN1UH1}
\begin{array} {ll}
G \ & =  \O(n)\cdot K_1N_1^+K_1\cdot U\cdot C\\
\   & = \O(n)\cdot N_1^+U\cdot K_1C\\
\   & = \O(n)\Sigma \cdot N_1^+U\cdot K_1C\\
\   & = \O(n)\cdot N_1^+U\cdot\Sigma K_1C\\
\   & = \O(n)\cdot N_1^+U\cdot H_1. 
\end{array}
\end{equation}

\subsection{Choice of Haar measures on subgroups of  \(G\).}  
\label{sec:Haar} 
Our next aim is to choose the Haar measures on each of the subgroups
defined in the previous section, so that the equalities
(\ref{eq:decompose}), (\ref{eq:G=KLUA}) and (\ref{eq:G=KN1UH1}) also
hold, in an appropriate sense, with respect to the products of the
chosen Haar measures.

\subsubsection*{Choice of Haar measure  \(\tilde \mu\) on
 \(G\).}

We choose a Haar integral $dk$ on $\O(n)$ such that $\vol(\SO(n))=1$;
in particular, 
\begin{equation} \label{eq:O(n):Haar}
\vol(\O(n)) = \int_{\O(n)} 1 \;dk = 2.
\end{equation}

We choose the Haar integral $d\bfn$ on $N$ (see~(\ref{eq:N}))such
that
\[
d\bfn=\prod_{i<j} dn_{ij}. 
\]

We choose the Haar integral $d\bfa$ on $A$ such that $\forall\, f\in\Cc(A)$,
\[
\int_A f(\bfa)\; d\bfa = 
\int_{(\R_{>0})^{n-1}} f(\bfa)\; \frac{da_1}{a_1}\cdots
\frac{da_{n-1}}{a_{n-1}}; \qquad (see
(\ref{eq:A})) 
\]
alternative notation: $d\bfa=\prod_{i=1}^{n-1} da_i/a_i$.
 
We choose a Haar measure $\tilde\mu$ on $G$ such that, 
\begin{equation} \label{eq:G:Haar}
\int_G f \; d\tilde\mu = \int_{\O(n)\times N\times A} f(k\bfn\bfa)\; dk\,
d\bfn\, d\bfa, \ \forall\,f\in\Cc(G).
\end{equation}

\subsubsection*{Decomposition of integrals on  \(A\) and
 \(N\).}  
We choose a Haar integral $d\bfc$ on $C$ such that (see~\ref{eq:C})
\[
d\bfc = (dc_1/c_1)\cdots(dc_{r_1+r_2-1}/c_{r_1+r_2-1}).
\]
Choose the Haar integral $d\bfa_1:=\prod_{i=1}^{r_2} d\beta_i/\beta_i$
on $A_1$ (see (\ref{eq:A1})). Then $d\bfa=d\bfa_1\, d\bfc$, where
$\bfa=\bfa_1\bfc$, $(\bfa_1,\bfc)\in A_1\times C$ (see
(\ref{eq:decompose})).

Let $d\bft$ denote the standard Lebesgue measure on $\R^{r_2}$.  Let
$\bfx$ denote the standard Lebesgue measure on $\cU$.  Then
$d\bfn=d\bft \, d\bfx$, where $\bfn=h(\bft) u(\bfx)$,
$(\bft,\bfx)\in \R^{r_2}\times \cU$.

\subsubsection*{Choice of Haar integral  \(dl\) on
 \(L_1\).}
 
Let $dl$ be a Haar integral on $L_1$ such that,
\begin{equation} \label{eq:L1:Haar}
\int_{L_1} f(l) \; dl = \int_{K_1\times \R^{r_2} \times A_1}
f(kh(\bft)\bfa_1)\; d\theta(k)\, d\bft \, d\bfa_1, \ \forall\,f\in\Cc(L_1),
\end{equation} 
where $\theta$ denotes the Haar measure on $K_1$ such that
\begin{equation} \label{eq:K1:Haar}
\theta(K_1)=1.
\end{equation}

\subsubsection*{Decomposition of Haar integral  \(d\tilde\mu\).}   
From the above choices of Haar integrals on various subgroups of $G$,
their interrelations, (\ref{eq:N(U)}) and (\ref{eq:G=KLUA}) we have
\begin{equation} \label{eq:G=KLUA:Haar} 
\int_G f(g)\; d\tilde\mu(g) = 
\int_{\O(n)\times L_1\times \cU\times C} 
f(kl\bfx\bfc)\; dk\, dl\, d\bfx\, d\bfc, \ \forall\, f\in\Cc(G).
\end{equation}

\subsubsection*{Choice of Haar measure  \(\tilde \nu\) on
 \(H\).}  
We also choose a Haar measure $\tilde\nu$ on $H$ such that for the
Haar measure $\tilde\nu_1:=\psi_\ast(\tilde \nu)$ on $H_1$ (see
(\ref{eq:tilde:nu1:def})), we have
\begin{equation} \label{eq:H1:Haar}
\int_{H_1} f\; d\tilde{\nu}_1 = 
\sum_{\sigma\in\Sigma} \int_{K_1\times C} f(\sigma k \bfc)\;
d\theta(k)\, d\bfc,\qquad \forall\, f\in\Cc(H_1).
\end{equation}

\subsection{Description of integral  \(\eta_1\) on
 \(G/H_1\).}
  \label{sec:eta1}
In order to describe $\eta_1$, we will express the integral
$d\tilde\mu$ as a product of an integrals on certain subset of $G$ and
the integral $d\tilde\nu_1$ using the expressions
(\ref{eq:G=KLUA:Haar}) and (\ref{eq:H1:Haar}).

\subsubsection*{A new description of the integral \(dl\).}  
First we will express the Haar integral on $L_1$ in terms of the product
decomposition $L_1=K_1N_1K_1$. 

By Proposition~\ref{prop:KNA=KNK} (stated and proved in
Appendix~\ref{append:KNA=KAK=KNK}), the following holds: $\forall\,
f\in\Cc(\SL_2(\R))$,
\begin{equation} \label{eq:KNA=KNK}
\begin{array}{l}
\int_{\SO(2)\times \R\times \R_{>0}} 
f(kh(t)\diag(\beta,\beta\inv))\; d\vartheta(k)\,dt\,(d\beta/\beta)\\[5pt] 
 \qquad = (\pi/2) \int_{\SO(2)\times \R_+\times \SO(2)}
f(k_1h(t^{1/2})k_2)\; d\vartheta(k_1)\,dt\,\vartheta(k_1),
  \end{array}
\end{equation}
where $\vartheta$ is a probability Haar measure on $\SO(2)$.

Since $L_1\cong\SL_2(\R)^{r_2}$, by (\ref{eq:L1:Haar}) and
(\ref{eq:KNA=KNK}), $\forall\, f\in\Cc(L_1)$,
\begin{equation} \label{eq:L1=KNK:Haar}
\int_{L_1} f(l)\; dl = (\pi/2)^{r_2}\int_{K_1\times (\R_+)^{r_2}\times K_1}
f(kh(\bft^{1/2})k')\; d\theta(k)\, d\bft\, d\theta(k'),
\end{equation}
where the notation is  
\begin{equation} \label{eq:bft:half}
\mbox{$\bft^{1/2}:=(t_1^{1/2},\ldots,t_{r_2}^{1/2})$, 
$\forall\,\bft=(t_1,\ldots,t_{r_2})\in(\R_+)^{r_2}$.}
\end{equation}

From (\ref{eq:G=KN1UH1}) and
(\ref{eq:G=KLUA:Haar})--(\ref{eq:L1=KNK:Haar}), $\forall\,
f\in\Cc(G)$,
\begin{eqnarray*} 
 &\ & \int_G f(g)\; d\tilde\mu(g)\\
 &=& (\pi/2)^{r_2}
\int_{\O(n)\times K_1\times (\R_+)^{r_2}\times K_1\times\cU\times C}
f(kk'_1h(\bft^{1/2})k_1u(\bfx)\bfc)\; {\scriptstyle \times} \\[3pt]  
 &\ & \qquad\ \qquad \ \qquad \qquad \ \qquad  {\scriptstyle\times}\;
dk\,d\theta(k'_1)\,d\bft\, d\theta(k_1)\, d\bfx\, d\bfc
\\
 &=& 
(\pi/2)^{r_2}(\#\Sigma)\inv \sum_{\sigma\in\Sigma}
\int_{\O(n) \times (\R_+)^{r_2}\times \cU \times K_1\times C} 
f(k\sigma h(\bft^{1/2})u(\bfx) k_1 \bfc)\;{\scriptstyle \times} \\
 &\ & \qquad\ \qquad \ \qquad \qquad \ \qquad {\scriptstyle \times}\; 
 dk\, d\bft\, d\bfx\, d\theta(k_1)\, d\bfc.
\\
 &=& \pi^{r_2}2^{-r_1-r_2}\int_{\O(n)\times(\R_+)^{r_2}\times
\cU\times H_1} f(k h(\bft^{1/2})u(\bfx) h_1)\; dk\, d\bft\, d\bfx\,
d\tilde\nu_1(h_1). 
\end{eqnarray*}
Now in view of (\ref{eq:eta1:def}), we have the following:

\begin{prop} \label{prop:G/H1:Haar}
For any $\bar f\in\Cc(G/H_1)$,
\[
\int_{G/H_1} \bar f\; d\eta_1
= (2\pi)^{r_2}2^{-n} \int_{\O(n)\times (\R_+)^{r_2}\times \cU} 
\bar f(kh(\bft^{1/2})u(\bfx) H_1)\; dk\, d\bft\, d\bfx.
\]
\end{prop}

\subsection{Changing the order of Integration.} \label{sec:R1} 

The Euclidean norm on $\M_n(\R)$ is invariant under the left and the
right multiplication by the elements of $\O(n)$. Therefore
\[
R^1_T=\O(n) \Psi(D^1_T) H_1/H_1,
\] 
where 
\begin{eqnarray}
\Psi(\bft,\bfx) &=& h(\bft^{1/2})\bfu(\bfx), \ \forall\,
(\bft,\bfx)\in(\R_+)^{r_2} \times \cU, \
\mbox{(see~(\ref{eq:bft:half}))} \nonumber \\
D^1_T &=& \{(\bft,\bfx)\in(\R_+)^{r_2}\times \cU: 
\|\con{\Psi(\sbft,\sbfx)}X_1\|<T\}
\label{eq:D1}
\end{eqnarray}

Since $\cU\cong \R^{\half n(n-1)-r_2}$, let $\ell$ denote the standard
Lebesgue measure on $(\R_+)^{r_2}\times \cU$. Then by
(\ref{eq:O(n):Haar}) and Proposition~\ref{prop:G/H1:Haar},
\begin{equation} \label{eq:R1:D1:vol}
\eta_1(R_T^1)=(2\pi)^{r_2}2^{-(n-1)}\ell(D^1_T).
\end{equation}

For the purpose of analysing the limit in Theorem~\ref{thm:F_T:weak},
we change the order of integration in Proposition~\ref{prop:R1} as follows:

\begin{prop} \label{prop:fubini2}
For all $f \in\Cc(G/\Gamma_1)$,
\[
\begin{array}{l} 
\frac{1}{\eta_1(R^1_T)}
 \int_{R^1_T} \left(\int_{H_1/H_1\cap\Gamma_1}
f(gh\Gamma_1)\; d\nu_1(\dot{h})\right) \,d\eta_1(\dot{g})
\\[7pt]
\qquad \ \qquad = 
(1/2)\int_{\O(n)} dk \cdot \int_{H_1/H_1\cap\Gamma_1}
 d\nu_1(\dot{h})\;
{\scriptstyle \times} \\[3pt]
\qquad \ \qquad \ \qquad {\scriptstyle \times}\;  
\left(\frac{1}{\ell(D^1_T)} \int_{(\sbft,\sbfx)\in D^1_T} 
f(k\Psi(\bft,\bfx)\Gamma_1)\; d\bft\, d\bfx \right). 
\end{array}
\]
\end{prop}

\subsection{Description of the set  \(D^1_T\).} 
\label{sec:D1} 

Our aim of this subsection is to show that $D^1_T$ is asymptotically
 the image of a ball of radius $T$ under a `polynomial like' map (see
 Propositions \ref{prop:D:B:asymptotic} and \ref{prop:Theta}).

\subsubsection*{Coordinates of \( \con{\Psi(\sbft,\sbfx)X_1} \).}    

Take $\bfx=(x_{ij})\in\cU$. If $u(\bfx)\inv=u(\bfy)$,
$\bfy=(y_{ij})\in\cU$, then
\[
y_{ij}= -x_{ij} + B_{ij}((x_{kl})_{0<l-k<j-i})
\]
where $B_{ij}:\prod_{0<l-k<j-i} M_{kl} \to M_{ij}$ is a polynomial
map for $i<j-1$, and $B_{ij}\equiv 0$ if $i=j-1$.
  
If $\con{u(\sbfx)}X_1=u(\bfx)X_1u(\bfy)=(\omega_{ij})_{i,j=1}^{r_1+r_2}$,
then $w_{ij}=0$ if $i>j$, and 
\[
\omega_{ij}=\left\{ 
\begin{array}{ll} 
   d_i & \mbox{if $i=j$ \ (see (\ref{eq:di}))}\\ 
   S_{ij}(x_{ij}) + Q_{ij}((x_{kl})_{0<l-k<j-i}) & \mbox{if $i<j$}, 
\end{array}  
\right.
\]
where $S_{ij}:M_{ij}\to M_{ij}$ ($i<j$) is defined as 
\begin{equation} \label{eq:Sij}
S_{ij}(x)=xd_j-d_ix, \ \forall\, x\in M_{ij},
\end{equation}
and $Q_{ij}:\prod_{0<l-k<j-i} M_{kl}\to M_{ij}$ is a polynomial map
for $i<j-1$, and $Q_{ij}\equiv 0$ if $i=j-1$.

Let $\bft=(t_i)\in(\R_+)^{r_2}$. If we write
\[
\con{h(\sbft)}\left(\con{u(\sbfx)}X_1\right)=(\zeta_{ij})_{i,j=1}^{r_1+r_2},
\] 
then $\zeta_{ij}=0$ if $i>j$, and
\[
\zeta_{ij}=
h(t_{i-r_1}^{1/2})\omega_{ij}h(-t_{j-r_1}^{1/2}) \qquad \mbox{if $i\leq j$}.
\]
where the convention is: $h(t_{i-r_1}^{1/2})=h(-t_{i-r_1}^{1/2})=1$ for
$1\leq i\leq r_1$.

Note that for $i=1,\ldots,r_2$, (see (\ref{eq:di}))
\[
h(t^{1/2})d_{r_1+i}h(-t^{1/2}) =\left( 
\begin{matrix} 
a_i-t^{1/2} b_i &  -(1+t) b_i
\cr 
b_i  & a_i +t^{1/2} b_i
\end{matrix}\right).
\]
Therefore
\[
\|\con{\Psi(\sbft,\sbfx)}X_1\|^2 = \|X_1\|^2 + 
\sum_{i=1}^{r_2} b_i^2(t_i^2+4t_i) + \sum_{i<j} |\zeta_{ij}|^2.
\]

\subsubsection*{Expressing \(D^1_T\) as an image of a ball.}
Now, in view of (\ref{eq:D1}), we want to find a function
\[
\tilde\delta: (\R_+)^{r_2} \times \cU \to (\R_+)^{r_2}\times\cU
\]
such that
\begin{equation} \label{eq:D:B}
\tilde\delta(B^+_{\scriptscriptstyle{\sqrt{T^2-\|X_1\|^2}}}\,)=D^1_T,
\end{equation}
where
\[
B^+_T:=\{(\bfs,\bfz)\in (\R_+)^{r_2}\times \cU:
\|\bfs\|^2+\|\bfz\|^2<T^2\}.
\]

Now for $(\bfs,\bfz)\in (\R_+)^r_2\times\cU$, we write
$\tilde\delta(\bfs,\bfz)=(\bft,\bfx)$, where
$\bft=(t_i)\in(\R_+)^{r_2}$ and $\bfx=(x_{ij})\in\cU$. Then
(\ref{eq:D:B}) holds, if we have:
\begin{eqnarray} 
\label{eq:si}
s_i &=& \sqrt{b_i^2(t_i^2+4t_i)}, \qquad (1\leq i\leq r_2)\\
\label{eq:zi}
z_{ij} &=& \zeta_{ij}, \qquad (1\leq i<j \leq r_1+r_2).
\end{eqnarray}

By first solving the equation (\ref{eq:si}), we get
\[
t_i=\sqrt{b_i^{-2}s_i^2+4}-2
\]
After that we solve the equation (\ref{eq:zi}) in the following order:
it is solved for all $\{(k,l):0<l-k<j-i\}$ before solving it for
the $(i,j)$. We get 
\begin{eqnarray}
\label{eq:xij}
x_{ij}&=&x_{ij}(\bft,\{x_{kl}:0<l-k<j-i\})\\
\ &=& S_{ij}\inv \left( 
             h(-t_{i-r_1}^{1/2})z_{ij}h(t_{j-r_1}^{1/2})
             - Q_{ij}((x_{kl})_{0<l-k<j-i})
                  \right). \nonumber
\end{eqnarray}

\subsubsection*{`Polynomial like' approximation for \(\tilde\delta\).} 

We put $\bft':=(t_i')\in(\R_+)^r_2$, 
\[
t_i'=|b_i|\inv s_i, \qquad 1\leq i\leq r_2.
\]
Next we put $\bfx':=(x_{ij}')\in\cU$, where (see~(\ref{eq:xij})) 
\[
x_{ij}'=x_{ij}(\bft',\{x_{kl}':0<k-l<j-i\}), \qquad (1\leq i<j\leq r_1+r_2). 
\]
Then we define 
\[
\delta(\bfs,\bfz)=(\bft',\bfx').
\]
It is straightforward to verify that   
\[
0\leq t_i'-t_i<2, \qquad 1\leq i\leq r_2.
\]
Therefore 
\begin{equation} \label{eq:approximate}
\delta(B_{T-2})\subset \tilde{\delta}(B_T) \subset \delta(B_T), \ 
\forall\, T>0.
\end{equation}

Also note that if $T>\|X_1\|$, then 
\[
T - \|X_1\|^2\, T\inv < \sqrt{T^2-\|X_1\|^2}\, < T.
\]
Therefore, since (\ref{eq:D:B}) and (\ref{eq:approximate}) hold, we
get the following:

\begin{prop} \label{prop:D:B:asymptotic}
For $T>\|X_1\|+2$, 
\[
\delta(B_{T-2-\|X_1\|^2T\inv}) \subset D_T^1 \subset
\delta(B_T).
\]
\end{prop}

\begin{prop} \label{prop:Theta}  
The map $\Theta:\R^{\half n(n-1)}\to G$ defined by
\[
\Theta(\bfs,\bfz):=\Psi(\delta((s_1^2,\ldots,s_{r_2}^2),\bfz)), \qquad 
\forall\,(\bfs,\bfz)\in \R^{r_2}\times\cU=\R^{\half{n(n-1)}}, 
\]
is a polynomial map; that is, each coordinate function of $\Theta$ is
a polynomial in $\half n(n-1)$-variables.
\end{prop}

\subsection{Jacobian of  \(\delta\).}
Let the notation be as in the definition of $\delta$. The Jacobian of
$\delta$ at $(\bfs,\bfz)$ is given by:
\begin{eqnarray} 
\Jac(\delta)(\bfs,\bfz)
  &:=& |\partial(\bft',\bfx')/\partial(\bfs,\bfz)| \nonumber \\
\  &=& \prod_{i=1}^{r_2} |\partial t_i'/\partial s_i| \cdot
      \prod_{i<j} |\partial x_{ij}'/\partial z_{ij}|
  \label{eq:lower:triang} \\
\  &=& \prod_{i=1}^{r_2} |b_i|\inv 
      \cdot \prod_{i<j} |\det(S_{ij})\inv|,
\label{eq:Jac(delta)}
\end{eqnarray}
where (\ref{eq:lower:triang}) holds because $\partial t_i'/\partial
z_{kl}=0$ for all $i,k,l$, and $\partial x_{kl}'/\partial z_{ij}=0$ for
all $0<l-k<j-i$, and (\ref{eq:Jac(delta)}) holds because $\det h(t)=1$
for all $t$. In particular, $\Jac(\delta)$ is a constant function. 

\subsubsection*{Computation of \(\det(S_{ij})\).}
By (\ref{eq:Mij:def})
\[
M_{ij}=\Hom(\R^{\nu_i},\R^{\nu_j})\cong \R^{\nu_i}\otimes
(\R^{\nu_j})^\ast,\qquad  (1\leq i<j\leq r_1+r_2), 
\]
where $\nu_k=1$ if $1\leq k\leq r_1$, and $\nu_k=2$ if $r_1< k \leq
r_2$. Under this canonical isomorphism, $S_{ij}$ corresponds to
\[
(1\otimes d_j^\ast)-(d_i\otimes 1),\qquad \mbox{(see (\ref{eq:Sij}))}
\]
whose eigenvalues are distinct, and by (\ref{eq:sigma:bar}) they are
\[
\sigma_{j'}(\alpha)-\sigma_{i'}(\alpha), \qquad i'\in \hat{i},\
j'\in\hat{j},
\]
where $\hat{k}=\{k\}$ if $\nu_k=1$, and $\hat{k}=\{k,r_2+k\}$ if
$\nu_k=2$. Therefore by (\ref{eq:Jac(delta)})
\begin{equation} \label{eq:J(P)}
\Jac(\delta) = 2^{r_2} \prod_{1\leq i<j\leq n}
  |\sigma_i(\alpha)-\sigma_j(\alpha)|\inv =
  2^{r_2}/\sqrt{|D_{\Q(\alpha)/\Q}|},
\end{equation}
where $D_{\Q(\alpha)/\Q}$ denotes the discriminant of $\Q(\alpha)$
over $\Q$.
\subsection{Volume of  \(R_T\).} \label{subsec:vol(R_T)}

We note that 
\begin{equation} \label{eq:vol(B_T)}
\ell(B^+_T)=2^{-r_2}\vol(B^{n(n-1)/2})T^{n(n-1)/2},
\end{equation}
where $\vol(B^m)$ denotes the volume of a unit ball in $\R^m$. Also
note that for any $m\in\N$ and $a,b>0$, if $T>\max\{a,b\}$ then
\[
((T+a)^m-(T-b)^m)/T^m<m(a+b)T\inv.
\]
Therefore by (\ref{eq:eta1:eta}), (\ref{eq:R1:D1:vol}),
Proposition~\ref{prop:D:B:asymptotic}, and since $\Jac(\delta)$ is a
constant,
\[
  \begin{array}{lcl}
\lim_{T\to\infty} \eta(R_T)/\ell(B^+_T)
&=& \lim_{T\to\infty} \eta_1(R^1_T)/\ell(B^+_T)\\
\ &=& (2\pi)^{r_2} 2^{-(n-1)} \lim_{T\to\infty} \ell(D_T^1)/\ell(B^+_T)\\
\ &=& (2\pi)^{r_2} 2^{-(n-1)} \Jac(\delta).
  \end{array}
\]
Now by (\ref{eq:J(P)}) and (\ref{eq:vol(B_T)}),
\begin{equation} \label{eq:R_T:vol}
c_\eta:=\lim_{T\to\infty} \eta(R_T)/T^{n(n-1)/2}=
\frac{(2\pi)^{r_2}\vol(B^{n(n-1)/2})}{2^{n-1}\sqrt{|D_{\Q(\alpha)/\Q}|}}.
\end{equation}

\section{Equidistribution of trajectories}

In view of Propositions~\ref{prop:R1} and \ref{prop:fubini2},
and since $\Jac(\delta)$ is a constant, for any
$f_1\in\Cc(G/\Gamma_1)$, and any $x_1\in G/\Gamma_1$,
\begin{equation} \label{eq:D1:B+:limit}
\begin{array}{l}
\lim_{T\to\infty}\frac{1}{\ell(D^1_T)} 
\int_{(\sbft,\sbfx)\in D^1_T} 
f_1(\Psi(\bft,\bfx)x_1)\; d\bft\,d\bfx \\[5pt]
\qquad = \lim_{T\to\infty}\frac{1}{\ell(B^+_T)} 
       \int_{B^+_T} f_1(\Theta(\bfs,\bfz)x_1)\; d\bfs\,d\bfz,
\end{array}
\end{equation}
where $\Theta$ as in Proposition~\ref{prop:Theta}. 

\begin{lema} \label{lema:dense}
For $x\in G/\Gamma_1$, if $H_1x$ is compact then $\cl{Ux}=G/\Gamma_1$.
\end{lema} 

\bp 
Choose $\bfc\in C$, such that $c_1>\ldots>c_{r_1+r_2}>0$
(see~\ref{eq:C}). Then $U=\{u\in G:\bfc^{-m} u \bfc^m\to 1 \mbox{ as
$m\to\infty$}\}$, which is the expanding horospherical subgroup of
$G^0$ associated to $\bfc$. Therefore by \cite[Prop.~1.5]{DR}
\begin{equation} \label{eq:dense}
\cl{\cup_{n=1}^\infty \bfc^m Uy}=G^0\Gamma_1/\Gamma_1=G/\Gamma_1,
\qquad \forall y\in G/\Gamma_1.
\end{equation}

Recall that $C\subset H_1$ and $H_1\subset N_G(U)$ (see
Section~\ref{sec:decompose}). Let $F$ be a compact subset of $H_1$
such that $Fx=H_1x$. Then by (\ref{eq:dense})
\begin{equation} \label{eq:FUx}
G/\Gamma_1=\cl{CUx}\subset
\cl{H_1Ux}=\cl{UH_1x}=\cl{UFx}=F\cl{Ux}.
\end{equation}

By Moore's ergodicity theorem~\cite{Moore}, $U$ acts ergodically on
$G/\Gamma_1$. Hence there exists $x_1\in G/\Gamma_1$ such that
$\cl{Ux_1}=G/\Gamma_1$. By (\ref{eq:FUx}) there exist $h\in F$ and
$x_2\in\cl{Ux}$ such that $x_1=hx_2$. Therefore, since $h\in N_G(U)$, 
\[
G/\Gamma_1=\cl{Ux_1}=\cl{Uhx_2}=h\cl{Ux_2}\subset h\cl{Ux}.
\]
Hence $\cl{Ux}=G/\Gamma_1$.
\ep

\begin{prop} \label{prop:polynomial}
For all $f_1\in\Cc(G/\Gamma_1)$, $k\in K$ and $h\in H_1$:
\[\lim_{T\to\infty} \frac{1}{\ell(B^+_T)} \int_{B^+_T}
       f_1(k\Theta(\bfs,\bfz)h\Gamma_1)\; d(\bfs,\bfz)
=\frac{1}{\mu_1(G/\Gamma_1)} \int_{G/\Gamma_1} f_1 \; d\mu_1,
\]
where $\Theta$ as in Proposition~\ref{prop:Theta}. 
\end{prop}

\bp Note that $G/\Gamma_1=G^0/(\Gamma_1\cap G^0)$ and $G^0=\SL_n(\R)$.
We apply Theorem~\ref{thm:polynomial} for $\Gamma_1\cap G^0$ in place
of $\Gamma$, $x=h\Gamma_1$ and the function $f_2\in\Cc(G/\Gamma_1)$,
where $f_2(g\Gamma_1):=f_1(kg\Gamma_1)$, $\forall\, g\in G$. Since
$H_1\Gamma_1/\Gamma_1=\bar\psi(H\Gamma/\Gamma)$, by
Proposition~\ref{prop:H:compact}, $H_1x$ is compact. Therefore by
Lemma~\ref{lema:dense}, $U_1x$ is dense in $G/\Gamma_1$. Since
$\Theta(\R^{r_2}\times\cU)\supset U$, the conclusion of
Theorem~\ref{thm:polynomial} holds, and hence the proposition follows.
\ep

\subsection{Proof of Theorem~\protect\ref{thm:EMS:counting}.} By a series of
reductions in Section~\ref{sec:2}, we showed that it is enough to
prove Theorem~\ref{thm:F_T:weak}. Now this result follows from
Propositions \ref{prop:R1} and \ref{prop:fubini2},
Equation~(\ref{eq:D1:B+:limit}), Proposition~\ref{prop:polynomial},
Lebesgue's dominated convergence theorem, Equation~(\ref{eq:nu1=nu}),
and the fact that $\mu_1=\bar{\psi}_\ast(\mu)$.  \ep

\section{Computation of  \(C_P\)}
\label{sec:volumes}

The rest of the article is devoted to proving the following:

\begin{theo} \label{thm:C_P}
Let the notation be as in Theorem~\ref{thm:EMS:counting}. Then
\[
C_P=\sum_{\cO\supset \Z[\alpha]}
\kappa(\cO)\cdot\frac{\vol(B^{n(n-1)/2})}{\vol(\SM_n)},
\]
where $\alpha$ is any root of $P$, the sum is over all orders $\cO$ of
the number field $K=\Q(\alpha)$ containing $\Z[\alpha]$,
\begin{eqnarray*}
\kappa(\cO)
&:=&\frac{2^{r_1}(2\pi)^{r_2}h_\cO R_\cO}{w_\cO\sqrt{|D_{K/\Q}|}},
\mbox{ here }\\
r_1   &=& \mbox{ Number of real places of $K$}, \\
r_2   &=& \mbox{ Number of complex places of $K$}, \\
h_\cO &=& \mbox{ Number of modules classes with order $\cO$} \\
R_\cO &=& \mbox{ Regulator of $\cO^\times$}, \mbox{ (see (\ref{eq:R_O}))} \\
w_\cO &=& \mbox{ Order of the group of roots of unity in $\cO^\times$},\\
D_{K/Q}   &=& \mbox{ Discriminant of $K$,} 
\end{eqnarray*}
(see \cite[pp.\ 10--17]{Koch} or Section~\ref{sec:Orders}) and 
\begin{eqnarray*}
\vol(B^m)&=& \pi^{m/2}/\Gamma(1+m/2)\\ 
\         &=& \mbox{ Volume of a unit ball in $R^m$ (we take $m=\half
           n(n-1)$),}  \\
\vol(\SM_n) &=& \prod_{s=2}^n \pi^{-s/2}\Gamma(s/2)\zeta(s)\\
\              &=& \mbox{ Volume of the determinant one surface} \\
\              &\ & \mbox{ in the  Minkowski fundamental domain $\cM_n$.}      \end{eqnarray*}
(see \cite[Sect.~4.4.4, Theorem~4]{Terras} or Section~\ref{sec:G/Gamma})
\end{theo}

The computation of $C_P$ depends on: (i) obtaining representatives, say
$X_0$, for each $\Gamma$-orbits in $V_P(\Z)$, and then (ii) computing
$\nu(H/H\cap\Gamma)$ for the $H$ and the $\nu$ associated to $X_0$, (iii)
computing $c_\eta$ (see (\ref{eq:vol(R_T)})), and also  (iv) computing
$\mu(G/\Gamma)$. We already know $c_\eta$ (see~\ref{eq:R_T:vol}). 

\subsection{Orbits under  \(\Gamma\) in  \(V_P(\Z)\).}
\label{sec:Gamma:orbits}

We now describe a result due to Latimer and MacDuffee~\cite{MacDuffee}
on a correspondence between the classes of matrices and classes of
ideals; here two matrices are said to be in the same equivalence
class if they are in the same $\Gamma$-orbit.

Fix any root $\alpha$ of $P$. Any (nonzero) ideal $I$ of $\Z[\alpha]$
is a free $\Z$-module of rank $n$. We say that ideals $I$ and $J$ of
$\Z[\alpha]$ are \emph{equivalent}\/ if and only if $aI= bJ$ for some
nonzero $a,b\in \Z[\alpha]$. Let $[I]$ denote the class of ideals in
$\Z[\alpha]$ equivalent to $I$.

For any $X\in V_P(\Z)$, $\alpha$ is an eigenvalue of $X$, and there
exists a nonzero eigenvector
$\bfomega:=\trn{(\omega_1,\ldots,\omega_n)}\in \Q(\alpha)^n$ such
that
\begin{equation} \label{eq:eigen}
X\bfomega=\alpha\bfomega
\end{equation}
Replacing $\bfomega$ by some integral multiple, we may assume that
$\omega_i\in \Z[\alpha]$ for $1\leq i\leq n$. Put
$I_{X}=\Z\omega_1+\cdots+\Z\omega_n$. Then by (\ref{eq:eigen}),
$\alpha I_{X}\subset I_{X}$. Hence $I_{X}$ is an ideal of
$\Z[\alpha]$. The ideal class $[I_{X}]$ depends only on $X$, and not
on the choice of the eigenvector $\bfomega$.

Now let $\gamma\in\Gamma$ and $Y=\con{\gamma}X$. Then
$\bfomega':=\gamma\bfomega\in I_{X}$, and
$Y\bfomega'=\alpha\bfomega'$. Let
$I_Y=\Z\omega'_1+\cdots+\Z\omega'_n$, where
$\trn{(\omega'_1,\ldots,\omega'_n)}:=\bfomega'$. Then $I_Y\subset
I_{X}$. Since $\gamma\inv\in \Gamma$, we have $\bfomega=\gamma\inv
\bfomega'\in I_Y$, and hence $I_{X}=I_Y$. Thus the ideal class
$[I_{X}]$ depends only on the $\Gamma$-orbit $\con{\Gamma}X$, and not on
the choice of its representative $X$.

\begin{theo} \label{thm:MacDuffee}
The assignment $\con{\Gamma}X\mapsto [I_{X}]$ is a one-to-one
correspondence between the collection of $\Gamma$-orbits in $V_P(\Z)$
and the collection of equivalence classes of ideals in $\Z[\alpha]$.
\end{theo} 

\subsubsection{Orders in  \(\Q(\alpha)\).} \label{sec:Orders}

A subring $\cO$ of the number field $K=\Q(\alpha)$ is called an {\emph
order}\/ in $K$, if its quotient field is $K$, $\cO\cap\Q=\Z$, and its
additive group is finitely generated. 

A free $\Z$-submodule of $K$ (additive) of rank $n=[\Q(\alpha):\Q]$ is
called a \emph{lattice}\/ in $K$; for example, any (nonzero) ideal of
$\Z[\alpha]$ is a lattice in $K$. Two lattices $\frakM$ and $\frakM'$ in $K$
are said to be \emph{equivalent}\/, if $a\frakM=b\frakM'$ for some nonzero
$a,b\in\Q(\alpha)$. Let $\bar{\frakM}$ denote the class of lattices
equivalent to $\frakM$. For ideals $I$ and $J$ of $\Z[\alpha]$, we have
$[I]=[J]\Leftrightarrow \bar I=\bar J$.

For a lattice $\frakM$ in $K$,
\begin{equation} \label{eq:Order:def} 
\cO(\frakM):=\{\beta\in K:\beta \frakM \subset \frakM\} 
\end{equation}
is an order in $K$, it is called the order of $\frakM$, and it depends
only on the class $\bar{\frakM}$.

Let $\cO$ be an order in $K$. Then by the class number theorem
\cite[Theorem~1.9]{Koch}, there are only finitely many 
classes of lattices in $K$ with order $\cO$. This number is called the
\emph{class number of $\cO$} and denoted by $h_\cO$.

The ring $\cO_K$ of algebraic integers in $K$ is an order. Any order
$\cO$ in $K$ is contained in $\cO_K$, and $[\cO_K:\cO]<\infty$.  Also
$\Z[\alpha]$ is an order in $K$, and hence there are only finitely
orders $\cO$ in $K$ with $\cO\supset \Z[\alpha]$.

\begin{prop} \label{prop:finite:orbits2}
The $\Gamma$-orbits in $V_P$ are in one-to-one correspondence with the
classes of lattices in $K$ whose orders contain $\Z[\alpha]$.

In particular, each order $\cO$ containing $\Z[\alpha]$ is associated
to $h_\cO$ distinct $\Gamma$-orbits in $V_P(\Z)$, and the number of
distinct $\Gamma$-orbits in $\cO$ equals $\sum_{\cO\supset\Z[\alpha]}
h_\cO$.
\end{prop}

\proof In view of Theorem~\ref{thm:MacDuffee}, to any $\Gamma$-orbit
$\con{\Gamma}_X$ in $V_P$, we associate the lattice class $\bar I_X$ of
an ideal $I_X$ in $\Z[\alpha]$. We associate $\bar{I_X}$ to the orbit
$\con{\Gamma}X$. We note that $\cO(I_X)\supset\Z[\alpha]$.

Conversely, let $\frakM$ be a lattice in $K$ such that $\cO(\frakM)\supset
\Z[\alpha]$. Then there exists a nonzero integer $a$ such that
$I:=a\frakM$ is an ideal of $\Z[\alpha]$. By Theorem~\ref{thm:MacDuffee},
there exists $X\in V_P$, such that $[I]=[I_X]$. Therefore
$\bar\frakM=\bar I_X$, and hence $\bar\frakM$ is associated to a unique
orbit $\con{\Gamma}X$, and $\cO(\frakM)=\cO(I_X)$. This proves the
one-to-one correspondence.

Now the second statement follows from the class number theorem for
orders.  \ep

\subsection{Compactness and volume of  \(H/(H\cap\Gamma)\).}
\label{sec:compact}

Fix $X_0\in V_P(\Z)$ and let the notation be as before. Put
\[
\zx=\{Y\in \M_n(\R): YX_0=X_0Y\}.
\]
Since $X_0\in \M_n(\Q)$, we have that $\zx$ is the real vector space
defined over $\Q$. That is, $\zx$ is the real span of
$\zx(\Q):=\zx\cap\M_n(\Q)$, and $\zx(\Q)\otimes_\Q \R=\zx$.

Let $\bfomega=\trn{(\omega_1,\ldots,\omega_n)}\in\Z[\alpha]^n$,
$\bfomega\neq 0$, be such that $X_0\bfomega=\alpha\bfomega$. Since all
eigenvalues of $X_0$ are distinct, there exists an $\R$-algebra
homomorphism $\lambda:\zx\to \C$ given by $Y\mapsto \lambda_Y$, such
that $Y\bfomega=\lambda_Y\bfomega$. Now if $Y\in \zx(\Q)$ then
$\lambda_Y\in \Q(\alpha)$.

Let $I_{X_0}=\Z\omega_1+\ldots+\Z\omega_n$. Then $I_{X_0}$ is an ideal
of $\Z[\alpha]$, and hence $I_{X_0}\otimes_\Z
\Q\cong\Q(\alpha)$. Therefore $\{\omega_1,\ldots,\omega_n\}$ are
linearly independent over $\Q$. Hence if $Y\in\zx(\Q)$ and
$Y\bfomega=0$, then $Y=0$. Thus
\[
\ker\,\lambda\cap\zx(\Q)=0.
\] 

Let $Y_\beta$ denote the matrix of the multiplication by
$\beta\in\Q(\alpha)$ on the $\Q$-vector space $I_{X_0}\otimes_{\Z}\Q$,
with respect to the basis $\{\omega_1,\ldots,\omega_n\}$. The map
$\beta\mapsto Y_\beta$ is a $\Q$-algebra homomorphism of $\Q(\alpha)$
into $\M_n(\Q)$. Since $Y_\alpha=X_0$, $Y_\beta\in\zx(\Q)$. Also
$\lambda_{Y_\beta}=\beta$. Hence $\lambda:\zx(\Q)\mapsto \Q(\alpha)$
is a $\Q$-algebra isomorphism. In particular,
\[
\zx(\Q)=\Q[X_0] \mbox{ and } \zx=\R[X_0].
\]

Note that for $Y\in\zx(\Q)$, $\lambda_YI_{X_0}\subset I_{X_0}
\Leftrightarrow Y\in\M_n(\Z)$. Therefore
\begin{equation} \label{eq:zx(Z)=O}
\zx(\Z):=\zx\cap\M_n(\Z)=\{Y\in\zx(\Q):\lambda_Y\in\cO(X_0)\},
\end{equation}
where $\cO(X_0)$ denotes the order of $I_{X_0}$ (see
(\ref{eq:Order:def})).

Recall the Notation~\ref{not:sigma}. Define $\sigma_i(\bfomega):=
\trn{(\sigma_i(\omega_1),\ldots,\sigma_i(\omega_n))}$. Then
$X_0\sigma_i(\bfomega)=\sigma_i(\alpha)\sigma_i(\bfomega)$. Let
\[
g_1=(\sigma_1(\bfomega),\ldots,\sigma_n(\bfomega))\in\M_n(\C).
\]
Then 
\[
g_1\inv X_0 g_1 = \diag(\sigma_1(\alpha),\ldots,\sigma_n(\alpha)),
\]
and all the entries of this diagonal matrix are distinct. Therefore
$g_1\inv \zx g_1$ is a diagonal matrix. We define functions $D_i$ on
$\zx$ by 
\[
g_1\inv Y g_1=\diag(D_1(Y),\ldots,D_n(Y)).
\]
Since $\zx=\R[X_0]$, and the $D_i$'s are $\R$-algebra homomorphisms,
we have $D_i(\zx)\subset \R$ for $1\leq i\leq r_1$, and by
(\ref{eq:sigma:bar}),
\[
D_{r_1+r_2+i}(Y)=\bar D_{r_1+i}(Y),\qquad (1\leq i\leq r_2).
\]
Therefore 
\begin{equation} \label{eq:Di:det} 
\det(Y)=\prod_{i=1}^n
|D_i(Y)|=\prod_{i=1}^{r_1+r_2}|D_i(Y)|^{\nu_i}, \qquad \forall\, Y\in\zx,
\end{equation}
where $\nu_k=1$ if $1\leq k\leq r_1$, and $\nu_k=2$ if $r_1< k \leq
i_2$.  Since $D_i(Y)=\sigma(\lambda_Y)$, $\forall\, Y\in\zx(\Q)$, we have
\[
\det(Y)=N_{\Q(\alpha)/\Q}(\lambda_Y), \ \forall\, Y\in\zx(\Q).
\]

Therefore by (\ref{eq:zx(Z)=O})
\begin{eqnarray} 
H &=& \{Y\in\zx:|\det(Y)|=1\} \label{eq:H:det1} \\
H(\Z)&=& H\cap \zx(\Z) \nonumber \\
\    &=& \{Y\in\zx(\Q): \lambda_Y\in
         |N_{\Q(\alpha)/\Q}(\lambda_Y)|=1,\,\cO(X_0)\} \nonumber\\
\    &=& \{Y\in\zx(\Q):\lambda_Y\in \cO(X_0)^\times\}
         \nonumber\\
\    &\cong& \cO(X_0)^\times; \nonumber
\end{eqnarray}
here $\cO(X_0)^\times$ denotes the multiplicative group of the order
$\cO(X_0)$ which is same as the multiplicative group of unit norm
elements in $\cO(X_0)^\times$.

\subsubsection{Dirichlet's Unit theorem and Compactness of
 \(H/H(\Z)\).} 

\begin{theo} \label{thm:H/H(Z):compact}
$H/H(\Z)$ is compact.
\end{theo}

\bp
Define $l:H\to\R^{r_1+r_2}$ as
\[
l(h)=(\nu_1\log|D_1(h)|,\ldots,\nu_{r_1+r_2}\log|D_{r_1+r_2}(h)|), \
\forall\, h\in H,
\]
where $\nu_i=1$ if $i\leq r_1$, and $\nu_i=2$ if $i>r_1$.

Let
\[
E=\{(x_1,\ldots,x_{r_1+r_2})\in\R^{r_1+r_2}:x_1+\cdot+x_{r_1+r_2} =0\}.
\]
Then, by (\ref{eq:Di:det}) and (\ref{eq:H:det1}), $l:H\to E$ is a
surjective homomorphism.

By (\ref{eq:decompose}) $H_1=\Sigma\cdot K_1\cdot C$ is a direct
product decomposition; let $p:H_1\to C$ denote the associated
projection. We define $l_1:C\to E$ by
\[
l_1(\bfc)=(\log c_1,\ldots,\log c_{r_1+r_2}), \qquad \mbox{(see (\ref{eq:C}))}
\]
and extend it to $H_1$ by $l_1(h)=l_1(p(h))$, $\forall\, h\in H_1$.

We note that $l_1(g_0 h g_0\inv)=l(h)$ for all $h\in H$. Therefore
\[ 
\ker\, l = g_0\inv (\ker\, l_1) g_0 = g_0\inv \Sigma K_1 g_0.
\] 
Hence $\ker(l)$ is compact. 

We define $\ell:\cO(X_0)^\times\to E$, by 
\begin{equation} \label{eq:ell:def} 
\ell(\lambda) = 
(\nu_1\log|\sigma_1(\lambda)|, \ldots,
\nu_{r_1+r_2}\log|\sigma_{r_1+r_2}(\lambda)|),\   
\forall\, \lambda\in \cO(X_0)^\times.
\end{equation}

Clearly, $l(Y)=\ell(\lambda_Y)$ for all $Y\in H(\Z)$. By Dirichlet
unit theorem \cite[Theorem~1.13]{Koch}, $\ell(\cO(X_0)^\times)$ is a
lattice in $E$. Therefore $l(H)/l(H(\Z))$ is compact. Since $\ker(l)$
is compact, this completes the proof.  \ep

\subsubsection{Computation of  \(\nu(H/H(\Z))\).}

Let $\pr:E\to R^{r_1+r_2-1}$ be the projection on the first
$r_1+r_2-1$ coordinate space. We choose a measure $m$ on $E$ such that
its image under $\pr$ is the standard Lebesgue measure on
$\R^{r_1+r_2}$. Let $\bar m$ denote the associated measure on
$E/\ell(\cO(X_0)^\times)$. We note that $l_1:C\to E$ preserves the
choices of the Haar integrals $d\bfc$ and $dm$. 

Let $\tilde{K}_1=\Sigma K_1$. In view of (\ref{eq:Sigma}) and
(\ref{eq:K1:Haar}), let $\tilde\theta$ be the Haar measure on $\tilde
K_1$ such that
\[
\tilde\theta(\tilde K)=\#(\Sigma)\theta(K_1)=2^{r_1}.
\]
Then by (\ref{eq:H1:Haar}), $q:\tilde K_1\backslash H\to C$, defined
as $\tilde K_1 h=p(h)$, is an isomorphism and it preserves the chosen
associated measures on both sides.

Therefore $l_1\circ q: \tilde K_1\backslash H_1 \to E$ is a group
isomorphism and preserves the chosen Haar measures on both sides. Note
that $H\cap\Gamma=H(\Z)$, and
\[
l_1(H_1\cap \Gamma_1)=l(H\cap\Gamma)=l(H(\Z))=\ell(\cO(X_0)^\times).
\]
Therefore we have an isomorphism, 
\[
\tilde K_1\backslash H_1/(H_1\cap\Gamma_1) \cong
 E/\ell(\cO(X_0)^\times)
\]
preserving the invariant measures on both sides. Now by
Theorem~\ref{thm:vol:K:G:Gamma} (stated and proved in
Appendix~\ref{append:K:G:Gamma}),
\begin{equation} \label{eq:m:nu1}
\nu_1(H_1/(H_1\cap\Gamma_1)) = \frac{\tilde\theta(\tilde
K_1)}{\#(\tilde K_1\cap(H_1\cap\Gamma_1))}\cdot \bar
m(E/\ell(\cO(X_0)^\times)).
\end{equation}

By the Dirichlet's unit theorem, let
$\{\epsilon_1,\ldots,\epsilon_{r_1+r_2-1}\}$ be a set of
generators of $\cO^\times$ modulo the group of roots of unity. Then
\[
\ell(\cO(X_0)^\times)=\oplus_{j=1}^{r_1+r_2-1} \Z\,\ell(\epsilon_j).
\]
Hence, by (\ref{eq:ell:def}),
\begin{equation} \label{eq:R_O}
\bar{m}(E/\ell(\cO(X_0)^\times)) =
|\det\left((\nu_i\log|\sigma_i(\epsilon_j)|)_{i,j=1}^{r_1+r_2-1}\right)| =:
R_{\cO(X_0)}, 
\end{equation}
which is called the \emph{regulator}\/ of the the order $\cO(X_0)$
(see~\cite[Sect.~1.3]{Koch}). 

We note that $g_0\inv(\tilde K_1\cap(H_1\cap\Gamma_1))g_0=\ker(l)\cap
H(\Z)\cong \ker(\ell)$, which is the group of roots of unity in
$\cO(X_0)$, and its order is denoted by
$w_{\cO(X_0)}$. Therefore,
\begin{equation} \label{eq:w_O}
\#(\tilde{K}_1\cap(H_1\cap\Gamma_1))=w_{\cO(X_0)}.
\end{equation}

Now from (\ref{eq:m:nu1})--(\ref{eq:w_O}) we obtain the following:

\begin{theo} \label{thm:vol:H:Gamma} 
Let $\cO(X_0)$ be the order of the ideal $I_{X_0}$ of $\Z[\alpha]$
which is associated to $X_0$ as in Theorem~\ref{thm:MacDuffee}. Then
\[
\nu(H/H\cap\Gamma)=\nu_1(H_1/H_1\cap\Gamma_1) = 
2^{r_1}R_{\cO(X_0)}/w_{\cO(X_0)}. 
\]
\end{theo}

\subsection{Volume of  \(G/\GL_n(\Z)\).} 
\label{sec:G/Gamma} 
The volume of $G/\GL_n(\Z)$ was computed by C.L.~Siegel. To use that
computation here we need to compare the Haar measure on $G$ chosen for
Siegel's computation with the one chosen in (\ref{eq:G:Haar}). Instead
it will be more convenient for us to use the volume computations as in
\cite[Section~4.4.4]{Terras}, which is also uses Siegel's formula.

\subsubsection*{The space \(\cP_n\) of positive \(n\times n\) matrices.} 
Let $\cP_n$ be the space of $n\times n$ real positive symmetric
matrices. Then $\GL_n(\R)$ acts transitively on $\cP_n$ by
\[
(g,Y)\mapsto \trn{g}Yg, \ \forall\, (g,Y)\in \GL_n(\R)\times\cP_n.
\]
We consider a $\GL_n(\R)$-invariant measure $\mu_n$ on $\cP_n$ defined
as follows: If we write $Y\in\cP_n$ as $Y=(y_{ij})$, $y_{ij}=y_{ji}$,
$y_{ij}\in\R$, then
\[
d\mu_n(Y)=|\det(Y)|^{-(n+1)/2}\prod_{i\leq j} dy_{ij}.
\]

Let $\SP_n=\{Y\in\cP_n: \det(Y)=1\}$. Then $G$ acts transitively on
$\SP_n$, and preserves the invariant integral $dW$ on $\SP_n$
which is defined as follows: If we write $Y\in \cP_n$ as $Y=t^{1/n}W$,
($t>0$, $W\in\SP_n$), then
\begin{equation} \label{eq:mu:dW}
d\mu_n(Y)=(dt/t)dW.
\end{equation}

\subsubsection*{Volume of Minkowski fundamental domain.}

Let $\SM_n$ denote the Minkowski fundamental domain for the action
of $\GL_n(\Z)$ on $\SP_n$. We have chosen $d\mu_n$, and $dW$ such
that by \cite[Section~4.4.4, Theorem~4, pp.~168]{Terras}, which uses
Siegel's method, we have the following:
\begin{equation} \label{eq:SMn:vol}
\vol(\SM_n):=\int_{\SM_n} 1\; dW = \prod_{k=2}^n
\pi^{-k/2}\Gamma(k/2)\zeta(k). 
\end{equation}

\subsubsection*{Comparing volume forms.}  Now we want to compare the
volume forms $dnda$ on $\O(n)\backslash G$ and $dW$ on $\SP_n$ with
respect to the map $\O(n)g\mapsto \trn{g}g$.

Put $D=\{\bfb=\diag(b_1,\ldots,b_n):b_i>0\}$. Choose the Haar integral
$d\bfb=\prod_{i=1}^n db_i/b_i$ on $D$. Then
\begin{equation} \label{eq:D:A}
\mbox{$d\bfb=(dt/t)\,d\bfa$, 
where $\bfb=t^{1/n}\/\bfa$, $t>0$, $\bfa\in A$.} 
\end{equation}
By direct computation of the Jacobian of the map 
\[
(\bfn,\bfb)\mapsto Y:=\trn{(\bfn\bfb)}(\bfn\bfb)
\]
from $N\times D\to \cP_n$, one has (\cite[Sec.4.1, Ex.24,pp.23]{Terras})
\begin{equation} \label{eq:mu:ND}
d\mu_n(Y)=2^n d\bfn\, d\bfb.
\end{equation}
By (\ref{eq:mu:dW}), (\ref{eq:D:A}) and (\ref{eq:mu:ND}), for
$\bfn\in N$ and $\bfa\in A$, we have  
\begin{equation} \label{eq:dW:NA}
\mbox{$dW=2^{n-1}d\bfn\, d\bfa$, where $W=\trn(\bfn\bfa)(\bfn\bfa)$.}
\end{equation}

If $d(\bar g)$ denotes the Haar integral on $\O(n)\backslash
G\cong AN$ associated to the Haar integrals $dg$ and $dk$, then by
(\ref{eq:G:Haar}),
\begin{equation} \label{eq:O(n):G:Haar}
\mbox{$d\bar g=d\bfn\,d\bfa$, where $\bar g=O(n)\bfn\bfa$,
$\bfn\in N$, $\bfa\in N$}.
\end{equation}

Now for any $f\in\Cc(\SP_n)$, by (\ref{eq:dW:NA}) and
(\ref{eq:O(n):G:Haar}), we have
\begin{equation}  \label{eq:W:dKg}
\int_{\SP_n} f(W)\; dW =
2^{n-1}\int_{\O(n)\backslash G} f(\trn{g}g)\;d{\bar g}.
\end{equation}

\subsubsection*{Relating \(\vol(\SM_n)\) and \(\vol(G/\GL_n(\Z))\).}
By (\ref{eq:W:dKg}), the map $\O(n)g\mapsto \trn{g}g$ from
$\O(n)\backslash G$ to $\SP_n$ is a right $G$-equivariant
diffeomorphism, and it preserves the invariant integrals $2^{n-1}
d{\bar g}$ and $dW$. We also note that $\O(n)\backslash G$ is
connected, and $Z(G)$ is the largest normal subgroup of $G$ contained
in $K$. Therefore by Theorem~\ref{thm:vol:K:G:Gamma} (stated and
proved in Appendix~\ref{append:K:G:Gamma}),
\[
2^{n-1}\mu(G/\GL_n(\Z))=\frac{\vol(\O(n))}{\#(Z(G)\cap
\GL_n(\Z))}\vol(\SM_n).
\]
By (\ref{eq:O(n):Haar}), $\vol(\O(n))=2$, and $\#(Z(G)\cap
\GL_n(\Z))=2$. Also $\Gamma=\GL_n(\Z)$. Thus by (\ref{eq:SMn:vol}), we
have the following:

\begin{theo} \label{thm:G:Gamma:vol}
\[
\mu(G/\Gamma)=2^{-(n-1)} \prod_{k=2}^n \pi(k/2)\Gamma(k/2)\zeta(k).
\]
\end{theo}

\subsection{Proof of Theorem~\protect\ref{thm:C_P}.}

By Proposition~\ref{prop:finite:orbits2}, there exists a finite set
$\calF\subset V_P(\Z)$, such that $V_P(\Z)$ is a disjoint union of the
orbits $\con{\Gamma} X_0$, $X_0\in \calF$. By Theorem~\ref{thm:EMS:orbit},
(\ref{eq:R_T}), and (\ref{eq:F_T}),
\[
C_P=\sum_{X_0\in\calF} C_{X_0}.
\]
By Theorem~\ref{thm:EMS:R_T},  
\[
C_{X_0}=c_\eta\cdot\frac{\mu(H/H\cap\Gamma)}{\nu(G\cap\Gamma)}.
\]
Let $\cO(X_0)$ denote the order in $\Q(\alpha)$ associated to the
$\Gamma$-orbit $\con{\Gamma}X_0$ as in 
Proposition~\ref{prop:finite:orbits2}. Then by (\ref{eq:R_T:vol}),
Theorem~\ref{thm:vol:H:Gamma}, and Theorem~\ref{thm:G:Gamma:vol},
\begin{eqnarray*}
C_{X_0} 
&=&
\frac{(2\pi)^{r_2}\vol(B^{n(n-1)/2})}{2^{n-1}\sqrt{D_{\Q(\alpha)/\Q}}}
\cdot \frac{2^{r_1}R_{\cO(X_0)}/w_{\cO(X_0)}}  
{2^{-(n-1)} \prod_{k=2}^n \pi^{-k/2}\Gamma(k/2)\zeta(k)}\\
\ &=& 
\frac{(2\pi)^{r_2}2^{r_1}
R_{\cO(X_0)}}{w_{\cO(X_0)}\sqrt{D_{\Q(\alpha)/\Q}}}  
\cdot  
\frac{\vol(B^{n(n-1)/2})}{\vol(\SM_n)}.
\end{eqnarray*}

This shows that $C_{X_0}$ depends only on $\cO(X_0)$. We recall that
$\cO(X_0)\supset \Z[\alpha]$. By
Proposition~\ref{prop:finite:orbits2}, for each order $\cO$ in $K$
containing $\Z[\alpha]$, there exist exactly $h_\cO$ number of $X_0\in
\calF$, such that $\cO(X_0)=\cO$. Therefore
\[
C_P=\sum_{\cO\supset\Z[\alpha]} 
\frac{(2\pi)^{r_2}2^{r_1} h_\cO R_{\cO}}
{w_{\cO}\sqrt{D_{\Q(\alpha)/\Q}}} \cdot 
\frac{\vol(B^{n(n-1)/2})}{\vol(\SM_n)}.
\]
\ep 

\subsubsection*{Proof of Theorem~\protect\ref{thm:Z[alpha]}.} 
By our hypothesis $\Z[\alpha]$ is the integral closure of $\Z$ in
$K=\Q(\alpha)$, and hence $\Z[\alpha]$ is the maximal order $\cO_K$ in
$K$. Now the theorem follows immediately from Theorem~\ref{thm:C_P}.
\ep
 
\appendix

\section{Decompositions of Haar integrals on  \(\SL_2(\R)\)} 
\label{append:KNA=KAK=KNK} 

Let
\begin{eqnarray*} 
h(t)&=&\left(\begin{matrix} 1 & t \cr 0 & 1 
\end{matrix}\right), \qquad \forall\, t\in\R \\
a(\lambda)&=& 
\left(\begin{matrix} \lambda & \ \cr \ & \lambda\inv
\end{matrix}\right), \qquad \lambda>0.\\
k(\theta) &=&
\left(
\begin{matrix}   
\cos(2\pi\theta) & -\sin(2\pi\theta) \cr \sin(2\pi\theta) &
\cos(2\pi\theta)
\end{matrix} 
\right), \qquad \theta\in \R/\Z.   
\end{eqnarray*}

First will compare the decompositions of Haar integrals on $\SL_2(\R)$
with respect to the Iwasawa decomposition and the Cartan decomposition.

\begin{prop} \label{prop:KNA=KAK}
For any $f\in \Cc(\SL_2(\R))$, 
\begin{eqnarray} \label{eq:KNA:KAK:main} 
\ &\ & 
\int_{(\R/\Z)\times \R\times \R_{>0}} 
f(k(\theta_1)h(t)a(\lambda))\;
d\theta_1\, dt\, \frac{d\lambda}{\lambda} \\
\ &=& (\pi/2) 
\int_{(\R/\Z) \times A \times (\R/\Z)} 
f(k(\theta_2)a(\alpha)k(\theta))\, |\alpha^2-\alpha^{-2}|\;
d\theta_2\, \frac{d\alpha}{\alpha}\, d\theta. \nonumber
\end{eqnarray}
\end{prop} 

\bp
Suppose
$g=k(\theta_1)h(t)a(\lambda)=k(\theta_2)a(\alpha)k(\theta)$. Then
\begin{equation} \label{eq:KNA:KAK:compare}
\trn{g}g=a(\lambda)\trn{h(t)}h(t)a(\lambda) = 
k(-\theta)a(\alpha^2)k(\theta). 
\end{equation}
Substituting $\beta:=\alpha^2$, $\mu:=\lambda^2$, and
$\phi=2\pi\theta$, from (\ref{eq:KNA:KAK:compare}) we get,
\begin{equation} \label{eq:mu:t:beta:phi} 
\begin{array}{l}
\mu=(1/2)(\beta+\beta\inv) + (1/2)(\beta-\beta\inv)\cos(2\phi)\\[3pt]
t=-(1/2)(\beta-\beta\inv)\sin(2\phi).
\end{array}
\end{equation}
Therefore 
\[
|\partial (\mu,t)/\partial(\beta,\phi)| =
 \frac{|\beta-\beta\inv|}{\beta}\mu.  
\]
Hence 
\begin{equation} \label{eq:KNA:KAK:Jacobian}
|\partial(\lambda,t)/\partial(\alpha,\theta)| = 
2\pi\frac{|\alpha^2-\alpha^{-2}|}{\alpha}\lambda.
\end{equation}

Then by (\ref{eq:KNA:KAK:compare}) and (\ref{eq:mu:t:beta:phi}) the
map
\begin{equation} \label{eq:KAK:KNA:map}
(\theta_2,\alpha,\theta)\mapsto (\theta_1,t,\lambda),
\end{equation}
is surjective if $0\leq \theta < 1/2$, and $\alpha\geq1$, and it is
injective if $0\leq \theta<1/2$ and $\alpha>1$. Therefore the map
(\ref{eq:KAK:KNA:map}) is a differentiable, surjective, its degree at
regular points is $4$, and its Jacobian is given by
(\ref{eq:KNA:KAK:Jacobian}). This gives (\ref{eq:KNA:KAK:main}).  \ep

\bigskip Next, we will show that
$\SL_2(\R)=\SO(2)h(\R)\SO(2)$, and express the Haar integral on on
$\SL_2(\R)$ with respect to this decomposition.

\begin{prop} \label{prop:KNK=KAK}
For any $f\in\Cc(\SL_2(\R))$, 
\begin{eqnarray} 
\label{eq:KNK:KAK:main}
\ &\ & \int_{\R/\Z\times \R_+\times \R/\Z} f(k(\phi')h(t)k(\phi))\;
d\phi'\, dt^2\, d\phi \\
\ & = & \int_{\R/\Z\times \R>0 \times \R/\Z} 
f(k(\theta') a(\alpha) k(\theta))\, |\alpha^2-\alpha^{-2}| \;
d\theta'\, \frac{d\alpha}{\alpha}\, d\theta.
\nonumber
\end{eqnarray}
\end{prop}

\bp
If we write $g=k(\phi')h(t)k(\phi)=k(\theta')a(\alpha)k(\theta)$, then
\begin{equation} \label{eq:KNK:KAK:compare}
\trn{g}g=k(\phi)\trn{h(t)}h(t)k(\phi)=k(\theta) a(\alpha^2) k(\theta).
\end{equation}
Therefore,
\begin{equation} \label{eq:KNK:KAK:trace}
\tr(\trn{g}g)=1+t^2=\alpha^2+\alpha^{-2}.
\end{equation}
Consider the change of variables $s:=t^2$, and $\beta:=\alpha^2$. Then
\[
\partial s = \frac{\beta-\beta\inv}{\beta} \partial\beta.
\]
Clearly, $\partial\phi/\partial\theta=1$, and
$\partial t/\partial \theta=0$. Therefore 
\[
|\partial(s,\phi)/\partial(\beta,\theta)|=\frac{|\beta-\beta\inv|}{\beta},
\]
and hence
\begin{equation} \label{eq:KNK:KAK:Jacobian}
|\partial(s,\phi)/\partial(\alpha,\theta)| = 
\frac{2|\alpha^2-\alpha^{-2}|}{\alpha}.
\end{equation}

By (\ref{eq:KNK:KAK:compare}) and (\ref{eq:KNK:KAK:trace}), we have
that the map
\[
(\theta',\alpha,\theta) \to (\phi',s,\phi)
\]
is surjective if $\alpha\geq 1$, and it is one-one if
$\alpha>1$. Therefore the map is a differentiable, surjective, its
degree at regular points is $2$, and its Jacobian is given by
(\ref{eq:KNK:KAK:Jacobian}). This gives (\ref{eq:KNK:KAK:main}).  \ep

\bigskip
From Proposition~\ref{prop:KNA=KAK} and
Proposition~\ref{prop:KNK=KAK}, we obtain the following:

\begin{prop} \label{prop:KNA=KNK}
For any $f\in\Cc(\SL_2(\R))$, 
\begin{eqnarray*}
\ &\ & \int_{\R/\Z\times \R \times \R_{>0}} f(k(\theta)h(s)a(\lambda)) \;
d\theta\, ds\, \frac{d\lambda}{\lambda} \\
\ & = &
(\pi/2)\int_{\R/\Z \times \R_+ \times \R/\Z} f(k(\phi')h(t)k(\phi)) \;
d\phi'\, dt^2 \, d\phi.
\end{eqnarray*}

\end{prop}

\section{A Lemma on volume of two sided quotients}
\label{append:K:G:Gamma}
Let $G$ be a Lie group and $\Gamma$ a lattice in $G$. Assume that we
are given a Haar measure on $G$, and we want to find the volume of
$G/\Gamma$. In many cases one can find a compact subgroup $K$ of $G$
such that $E=K\backslash G$ is diffeomorphic to a Euclidean space, and
construct a fundamental domain, say $\cF$, for the right
$\Gamma$-action on $E$. The following result expresses the volume of
$G/\Gamma$ in terms of the volume of $\cF$.

\begin{theo} \label{thm:vol:K:G:Gamma}
Let $G$ be a Lie group and $K$ be a compact subgroup of $G$ such that
$K\backslash G$ is connected.  Let $\Gamma$ be a discrete subgroup of
$G$. Let $\tilde\mu$ (resp.\ $\nu$) be a Haar measures on $G$ (resp.\
$K$).  Let $\eta$ (resp.\ $\mu$) be the corresponding $G$-invariant
measure on $K\backslash G$ (resp.\ $G/\Gamma$). Let $\cF$ be a
measurable fundamental domain for the right $\Gamma$-action on
$K\backslash G$; in other words, $\cF$ is measurable and it is the
image of a measurable section of the canonical quotient map
$K\backslash G\to K\backslash G/\Gamma$. Then
\begin{equation} \label{eq:thm:K:G:Gamma}
\mu(G/\Gamma)=\frac{\nu(K)}{\#(K_0\cap\Gamma)}
\cdot\eta(\cF), 
\end{equation}
where $K_0$ is the largest normal subgroup of $G$ contained in $K$.
\end{theo} 

To prove this result, we need the the following two observations. 

\begin{lema} \label{lema:K_0}
For $\gamma \in G$, put 
\[
X_\gamma=\{\omega\in G: \omega \gamma \omega\inv \in K\}.
\]
Then either $X_\gamma$ is a finite union of strictly lower dimensional
analytic subvarieties of $G$, or $\gamma\in K_0$.
\end{lema}

\bp Because the map $\omega \mapsto \omega \gamma\omega\inv$ on $G$ is
an analytic map, and $K$ is a Lie subgroup of $G$, we have that
$X_\gamma$ is a finite union of analytic subvarieties of $G$. Therefore
either $X_\gamma$ is strictly lower dimensional, or $X_\gamma=G^0$. In
the latter case, since $KX_\gamma=X_\gamma$ and $KG^0=G$, we have
$X_\gamma=G$.

Put $K'=\{\gamma\in G: X_\gamma=G\}$. Then $K'$ is a normal
subgroup of $G$, and $K'\subset K$. Hence $K'\subset K_0$. This
completes the proof.  \ep
 
\begin{lema} \label{lema:K_0:Gamma}
Let $\Gamma$ be a discrete subgroup of $G$. Define 
\[ 
\mbox{
$K(g)=K\cap g\Gamma g\inv$ \  and \  $f(g)=\#(K(g))$, \ $\forall\, g\in G$.
}
\]
Then for $\tilde\mu$--a.e.\  $g\in G$, we have
\begin{equation} \label{eq:K_0:Gamma}
\mbox{$K(g)=g(K_0\cap\Gamma)g\inv$ \  and \  $f(g)=\#(K_0\cap\Gamma)$.}
\end{equation}
\end{lema}

\bp We put $n_0=\#(K_0\cap\Gamma)$. Since $K_0$ is normal in $G$ and
$K_0\subset K$, 
\begin{equation} \label{eq:K0:min}
K(g)\supset K_0\cap g\Gamma g\inv = g(K_0\cap\Gamma)g\inv,\
\forall\, g\in G.
\end{equation}

Take any $g\in G$. Since $K$ is compact and $\Gamma$ is discrete,
there exists an open neighbourhood $\Omega$ of $e$ in $G$ such that 
\[
\Omega K\Omega\inv \cap g\Gamma g\inv = K\cap g\Gamma g\inv.
\]
Therefore
\begin{equation} \label{eq:Omega}
K(\omega g)=\omega(\omega\inv K\omega \cap g\Gamma g\inv )\omega\inv 
\subset \omega K(g) \omega\inv, \ \forall\, \omega\in\Omega.
\end{equation}

First suppose, $f(g)\leq n_0$. Then by (\ref{eq:K0:min}) $n=n_0$, and
by (\ref{eq:Omega}),
\[
K(\omega g)=\omega K(g) \omega\inv=\omega g(K_0\cap\Gamma)g\inv
 \omega\inv,\  \forall\,\omega\in\Omega. 
\]
In particular, $f(\omega g)=n_0$ for all $\omega\in\Omega$.

Now suppose $f(g)>n_0$. Then by (\ref{eq:Omega})
\begin{eqnarray*}
\Omega g \cap f\inv(f(g)) &=& \{\omega g\in\Omega g: K(\omega g)=\omega
g(g\inv K g\cap\Gamma)g\inv \omega\inv\} \\
\ &\subset& \cap_{\gamma \in g\inv K g\, \cap\Gamma}\, X_\gamma.
\end{eqnarray*}
Now, by Lemma~\ref{lema:K_0}, either there exists $\gamma \in g\inv Kg
\cap\Gamma$ such that $X_\gamma$ is a finite union of strictly lower
dimensional analytic subvarieties of $G$, or $g\inv Kg\cap \Gamma
\subset K_0$. In the latter case, by (\ref{eq:K0:min}), $K(g)=
g(K_0\cap\Gamma)g\inv$, and hence $f(g)=n_0$, which is a
contradiction.

Thus we have shown that (i) for all $g\in f\inv(n_0)$, 
(\ref{eq:K_0:Gamma}) holds; and (ii) $\cup_{n\neq n_0} f\inv(n)$ is
contained in a countable union of strictly lower dimensional analytic
subvarieties of $G$, and hence $\tilde \mu(\cup_{n\neq n_0}
f\inv(n))=0$.  This completes the proof.  \ep

\subsubsection*{Proof of Theorem~\protect\ref{thm:vol:K:G:Gamma}.}
Consider the map $\psi:G/\Gamma\to K\backslash G/\Gamma$. For any $g\in G$
and $x=g\Gamma\in G/\Gamma$, we have
\[
\psi\inv(Kg\Gamma)=Kx\cong K/K\cap(g\Gamma g\inv)=K/K(g).
\]
Since $K(kg)=K(g)$, $\forall \, k\in K$, we can define $f(Kg):=f(g)$,
$\forall\, g\in G$. Now by Fubini's theorem,
\begin{equation} \label{eq:nu(K):f}
\mu(G/\Gamma)=\int_{Kg\in \cF} \nu(K)/f(Kg)\; d\eta(Kg).
\end{equation}

By Lemma~\ref{lema:K_0:Gamma}, $f(g)=\#(K_0\cap \Gamma)$ for $\tilde
\mu$--a.e.\ $g\in G$. Hence $f(Kg)=\#(K_0\cap\Gamma)$ for $\eta$--a.e.\
$Kg\in K\backslash G$. Now (\ref{eq:thm:K:G:Gamma}) follows from
(\ref{eq:nu(K):f}).  \ep

\subsubsection*{Acknowledgement.} 
The author would like to acknowledge that his discussions with G.\ A.\ 
Margulis, T.\ N.\  Venkataramana, S.\  Mozes and A.\  Eskin have
contributed significantly to this article.

\end{document}